\documentclass[12pt,a4paper]{article}

\usepackage[utf8]{inputenc}
\usepackage[T2A]{fontenc}
\usepackage[english]{babel}
\usepackage{amscd,amssymb}
\usepackage{amsfonts,amsmath,array}

\usepackage{longtable,wrapfig}
\usepackage{amsthm}
\usepackage{tikz}
\usepackage{hyperref}

\newtheorem{defin}{Definition}
\newtheorem{proposition}{Proposition}

\newtheorem{corollary}{Corollary}
\newtheorem{theorem}{Theorem}
\newtheorem{observation}{Observation}
\newtheorem{condition}{Condition}

\newtheorem{lemma}{Lemma}
\newtheorem{remark}{Remark}
\newtheorem{question}{Question}

\title{On the chromatic number of the plane for map-type colorings} 
\author{Georgy Sokolov, Vsevolod Voronov}

\begin{document}

\maketitle

\begin{abstract}
     We consider the Hadwiger--Nelson problem on the chromatic number of the plane under conditions of coloring a map containing a finite number of vertices in any bounded region. Woodall (1973) and Townsend (1981) showed that at least 6 colors would be required. In the present paper, it is shown that at least 7 colors are required to color a map in which the boundaries are not arcs of a unit circle and three boundaries connect at each vertex.  As a corollary, we obtain that at least 7 colors are required for a proper coloring in which the regions are arbitrary polygons. The proof relies on techniques developed for a similar result concerning the chromatic number of the plane with a forbidden interval of distances.
\end{abstract}

\maketitle

\section{Introduction}

In this paper, we consider one of the variations of the Hadwiger--Nelson problem on the chromatic number of the plane. The classical formulation of the problem requires finding the smallest number of colors that would be required to color the Euclidean plane $\mathbb{R}^2$ so that two points at unit distance have different colors. This number is called the chromatic number of the plane and is denoted by $\chi(\mathbb{R}^2)$. From \cite{deGrey,exoo2020chromatic} it is known that
\[
    5 \leq \chi(\mathbb{R}^2) \leq 7.
\]

There are several formulations of the problem in which tighter constraints on the coloring are imposed. In particular, one can ask: What is the number of colors required if the plane is partitioned into some monochromatic regions? In addition, one can define the kind of these regions. If we require that the sets of points colored in each color are measurable, then we have the problem of the measurable chromatic number of the plane, $\chi_{mes}(\mathbb{R}^2)$, for which Falconer  \cite{falconer1981realization}  proved the lower bound,
\[
    5 \leq \chi_{mes}(\mathbb{R}^2) \leq 7.    
\]
At this point, the best lower bound for the measurable chromatic number proved in 1981 is a consequence of the general result of \cite{deGrey}. 

Strengthening the constraints, we can require that the boundaries between regions are continuous images of segments (Jordan arcs). We denote such a chromatic number by $\chi_{map}(\mathbb{R}^2)$. In \cite{woodall1973distances,townsend2005colouring}, in addition, it is required that in any bounded region the number of vertices of the map is finite, and then it is possible to show that
\[
6 \leq \chi_{map}(\mathbb{R}^2) \leq 7. 
\]

Further, we can assume that the boundaries are piecewise--linear, the monochromatic regions are polygons, and, in addition, some additional conditions are satisfied, such as the restriction on the area of the monochromatic regions. In \cite{coulson2004chromatic}, a proof for this special case was independently given, namely, 
\[
    6 \leq \chi_{poly}(\mathbb{R}^2) \leq 7.    
\]

Trivially,
\[
    \chi_{\{1\}}(\mathbb{R}^2) \leq \chi_{mes}(\mathbb{R}^2) \leq \chi_{map}(\mathbb{R}^2) \leq \chi_{poly}(\mathbb{R}^2).
\]

In \cite{manta2021triangle}, it was shown that it would take 7 colors to color a plane divided into triangles.

In 2003, G. Exoo proposed a formulation of the chromatic number problem of a plane with a forbidden distance interval~\cite{exoo1}. It is required that there is no pair of monochromatic points $x,y$ such that the Euclidean distance $\|x-y\|$ belongs to the interval $[1 - \varepsilon, 1+\varepsilon]$. We denote the corresponding chromatic number $\chi_{[1 - \varepsilon, 1+\varepsilon)}(\mathbb{R}^2]$. In \cite{voronov2024chromatic}, it was shown that when $0<\varepsilon<1.15$ the equality is satisfied 

\[
    \chi_{[1 - \varepsilon, 1+\varepsilon]}(\mathbb{R}^2)=7.
\]

In this paper we extend this equality to the cases of $\chi_{poly}(\mathbb{R}^2)$ and $\chi_{map}(\mathbb{R}^2)$ under some additional conditions. In the first case, any bounded region must contain a finite number of polygons; in the second case, we require that any unit circle intersects the boundaries of the regions by a finite set of points, and in addition, we impose some restrictions on the vertices. The book \cite{Soifer} and the introduction to the article \cite{warsaw} are devoted to a review of other problems related to the chromatic number of the plane.

\section{Preliminaries}

\begin{condition}{(``proper'', ``Hadwiger--Nelson'')}.
    A proper $k$-coloring of the plane is a function $c: \mathbb{R}^2 \to \{ 1, 2, \dots, k \}$ such that for any points $x, y \in \mathbb{R}^2$ with $\|x - y\| = 1$, we have $c(x) \neq c(y)$.
    \label{cond_proper}
\end{condition}

In the following, we consider proper colorings of the plane throughout the paper. For the sake of brevity, we will sometimes omit the word ``proper''.

\begin{condition}{(``Jordan'')}

    \label{cond_jordan}
    A Jordan (proper) coloring of the plane is a coloring in which the plane is partitioned into vertices, boundaries, and regions, satisfying the following conditions:
\begin{itemize}
    \item  Vertices are points.
    \item Boundaries are Jordan curves terminating at vertices.
    \item Regions are the connected components formed by the vertices and boundaries, with all points in a given region assigned the same color.
    \item Any two regions sharing a boundary must be assigned different colors.
    \item Every vertex is an endpoint of at least three boundaries.
    \end{itemize}
\end{condition}

\begin{condition}{(lf = ``locally finite'')}
    A Jordan coloring of a plane is called locally finite if any bounded subset of the plane intersects only a finite number of boundaries.
    \label{cond_loc_finite}
\end{condition}

In the following, we consider only locally finite colorings. Therefore,  let us combine the previous two conditions.

\begin{condition}{(map = ``map-type coloring'')}
    A map-type coloring is a Jordan coloring that is also locally finite.
    \label{cond_map}
\end{condition}

\begin{defin}
    The degree of a vertex is the number of boundaries whose ends are a vertex.
\end{defin}

If we consider that the partition into regions was obtained as an embedding in the plane of some planar graph $G$, then the degree of a vertex $v$ is its degree as a vertex of $G$.

\begin{defin}
The multicolor of a point $u$, denoted $\overline{c}(u)$, is the set of colors assigned to the regions whose closure contains $u$.
\end{defin}

\begin{defin}
The chromaticity of a point $u$ is the number of elements in its multicolor $\overline{c}(u)$. A point of chromaticity $1$, $2$ and $3$ will be called monochromatic, bichromatic and trichromatic, respectively.
\end{defin}

\begin{observation}

    Without loss of generality, we assume that:
\begin{itemize}
    \item Each vertex is assigned the color of one of the regions whose closure contains it.
    \item Each boundary is assigned the color of one of the two regions it separates.
    \end{itemize}
\end{observation}

Note that all points inside the regions are monochromatic, points on the boundary of the regions are either bichromatic, or monochromatic (if the regions separated by the boundary are colored in the same color).
 In addition, the chromaticity of the vertex does not exceed its degree. The chromaticity may be less than the degree of a vertex if the colors of some regions whose closures contain the vertex coincide.

Next, we introduce a few specific properties that will be needed in the formulation of the main results.

\begin{condition}{(fa = ``forbidden arcs'')}
    If any boundary intersects with a circle of unit radius at a finite set of points, we say that arcs of unit curvature are forbidden.
    \label{cond_curv1}
\end{condition}

\begin{condition}{(3col = ``no trichromatic vertices of degree 4+'')}
    There are no trichromatic vertices of degree greater than 3.
     \label{cond_3col}
 \end{condition}

\begin{condition}{(cubic = ``no vertices of degree 4+'')}
    The degree of any vertex is 3, or, in other words, the map is an embedding of some planar cubic graph.
    \label{cond_cubic}
\end{condition}

\begin{condition}{(poly = ``polygonal'')}
    Polygonal coloring of the plane is such a map-type coloring that all regions are polygons.
    \label{cond_poly}
\end{condition}

\begin{observation} \label{l0.1}
    If arcs of unit curvature are forbidden, then the following holds:
    \begin{enumerate}
        \item The boundaries do not contain arcs of a circle of radius $1$.
        \item Any circle of radius $1$ is partitioned into a finite number of monochromatic arcs.
    \end{enumerate}
\end{observation}

In particular, these properties are satisfied for any polygonal coloring.

Let us list some notations used below.
\begin{itemize}
    \item Denote by $C_i \subset \mathbb{R}^2$ the set of points colored by $i$-th color. Sometimes for brevity we will use the notation $C_{123} = C_1 \cup C_2 \cup C_3$.

    \item The closure of some set $P \subset \mathbb{R}^2$ in the sense of standard topology will be denoted by $\overline{P}$.
    
    \item $\omega(x) \subset \mathbb{R}^2$ is a unit circle centered at the point $x$. 
    
    \item We will say that points $u$, $v$ are connected by a curve $uv$ if the context makes it clear which curve we are talking about, or if we are considering an arbitrary continuous curve with ends at points $u$, $v$.    

    \item Considering a neighborhood $U(x)$ of a point $x$, we will always take a neighborhood small enough to be contained in the union of regions whose color belongs to the multicolor of the point. 

\end{itemize}

\section{Main result}

\begin{theorem} \label{t2}

If the coloring is proper, Jordan, locally finite, if a vertex with degree greater than 3 cannot be trichromatic, and if arcs of unit curvature are forbidden, then at least 7 colors will be required, i.e. 
\[
\chi_{map+fa+3col}(\mathbb{R}^2) = 7.
\]
\end{theorem}

Obviously, Condition \ref{cond_cubic} is stronger than Condition \ref{cond_3col}. The formulation is somewhat simplified. 

\begin{corollary}
\[
\chi_{map+fa+cubic}(\mathbb{R}^2) = 7.
\] 
\end{corollary}

The conditions can be greatly simplified for a map whose regions are polygons. Some auxiliary steps are required in the proof. 

\begin{theorem} \label{t1}
There is no polygonal coloring of the plane in $6$ colors, i.e.
\[
\chi_{poly}(\mathbb{R}^2) = 7.
\]

\end{theorem}

\section{Simple properties of the coloring}

We assume that the plane is colored in $6$ colors and Conditions \ref{cond_proper}--\ref{cond_curv1} are satisfied.  Our goal is to demonstrate the impossibility of such a coloring under the conditions of Theorem \ref{t2}.

\begin{proposition} \label{l1.1}
Let $u$ and $v$ be two points of the same color, connected by a continuous monochromatic curve $uv$.
Then, there does not exist a point $z$ of the same color such that $\|u-z\|<1$, $\|v-z\|>1$.
\end{proposition}

\begin{proof} If there were such a point $z$, then $\omega(z) \cap uv \neq \emptyset$.
\end{proof}

\begin{proposition} \label{l1.2}
Let $u$ and $v$ lie on the boundary of two monochromatic regions.  Then, no point  $z$  satisfying $\|z-u\|<1$, $\|z-v\|>1$  can be colored with either of the two colors of these regions.
\end{proposition}

\begin{proof} Assume, for contradiction, that such a point $z$ exists. Without loss of generality, suppose that $u, v$ lie on the boundary of the regions of colors 1, 2. Then in arbitrarily small neighborhoods $u$, $v$ there are points $u_1$, $v_1$ of color 1 and points $u_2, v_2$ of color 2. We can choose these points such that $\||z-u_i\|<1$, $\|z-v_i\|>1$, $i=1,2$. It remains to apply the Proposition \ref{l1.1} to the pairs of points $u_1, v_1$ and $u_2, v_2$.\end{proof} 

\begin{proposition} \label{l1.3}
If a point $x$ is monochromatic, then for any point $y$ satisfying  $\|x-y\|=1$, the color of $x$ does not appear in the multicolor of $y$, i.e., $c(x) \in \overline{c}(y)$.
\end{proposition}

\begin{proof}
Suppose that there exists a monochromatic point $x$, $c(x)=1$, and a point $y$ such that $\|x-y\|=1$ and $1 \in c(y)$.
Let $z$ and $w$ be such points in $U(x)$ that $z$ lies at a distance less than $1$ from $y$, and $w$ lies at a distance greater than $1$. By definition of multicolor, in any neighborhood of $y$ there is a point of color $1$. Given a sufficiently small neighborhood, this point will also be at a distance less than $1$ from $z$ and greater than $1$ from $w$, which contradicts Proposition \ref{l1.1}. 
\end{proof}

\begin{proposition} \label{l1.4}  
Let \( x \) and \( y \) be points such that \( \|x - y\| = 1 \). Suppose \( x \) is connected by continuous curves of color 1 to points \( z \) and \( w \), where \( \|y - z\| > 1 \) and \( \|y - w\| < 1 \). The point \( x \) itself may have a different color. Then, color 1 does not appear in the multicolor of \( y \), i.e., \( 1 \notin \overline{c}(y) \).  
\end{proposition}  

\begin{proof}  
Assume, for contradiction, that \( y \) contains color 1 in its multicolor. Then, in any neighborhood of \( y \), there exists a point of color 1 that does not lie on the circle centered at \( x \) with radius 1.  

Let \( u \) be such a point in a sufficiently small neighborhood of \( y \), chosen so that   
 \( \|u - z\| > 1 \),  
 \( \|u - w\| < 1 \).  

Consider the unit circle $\omega(u)$. This circle must intersect the curve \( zxw \) at some point $x'$ other than \( x \), since by construction, \( \|u - x\| \neq 1 \). Observe that $c(u)=c(x')=1$, so the coloring is not proper.
\end{proof}

\begin{lemma} \label{l1.5}  
There do not exist two bichromatic points with the same multicolor that are exactly one unit apart.  
\end{lemma}  

\begin{figure}[ht!]
    \centering
    \includegraphics[width=12cm]{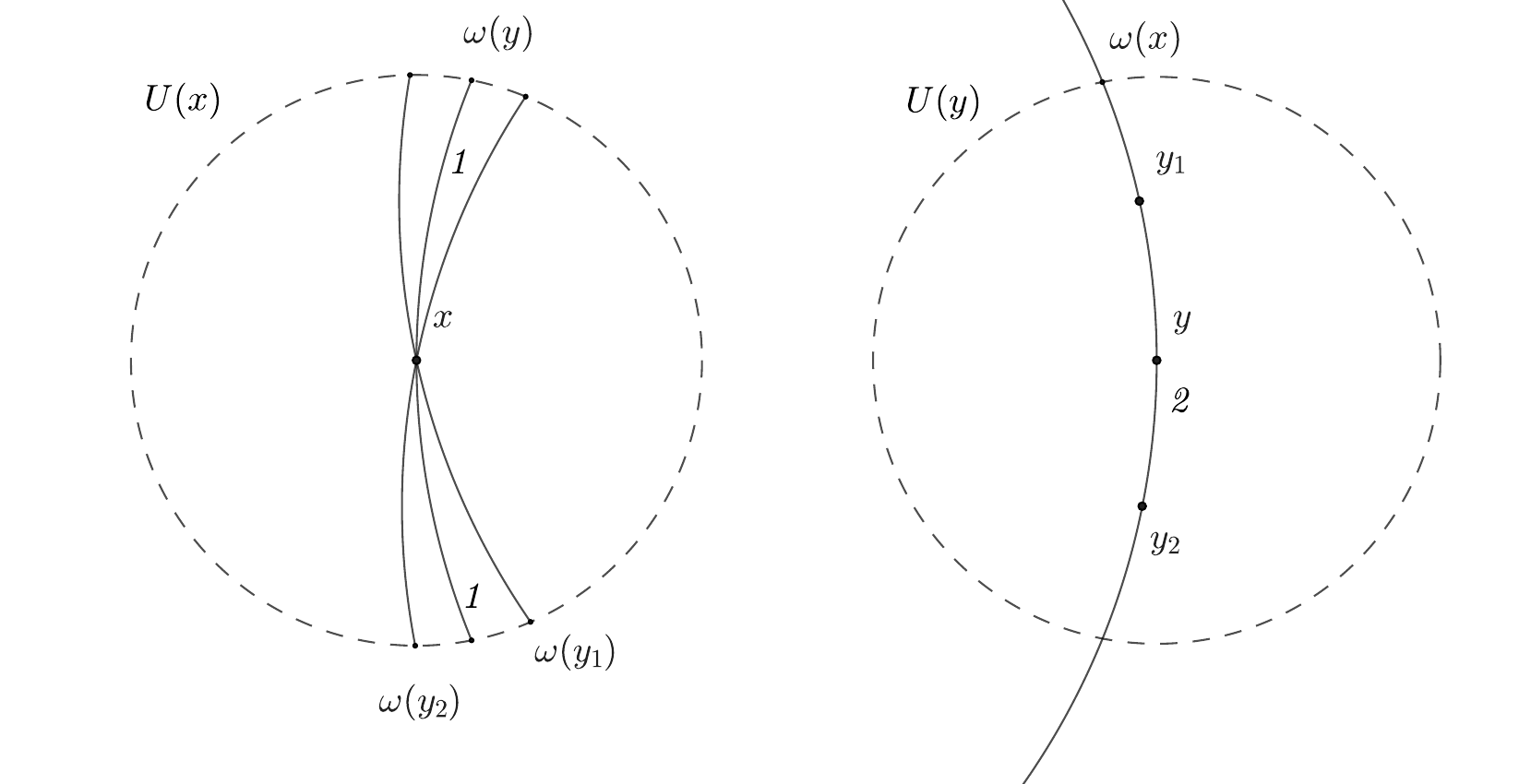}
    \caption{neighborhoods of two bichromatic points}
    \label{fig:two_2c_pts}
\end{figure}

\begin{proof}  
Assume, for contradiction, that such points exist. Let \( x \) and \( y \) be two bichromatic points at distance 1 apart, both having multicolor \( \{1,2\} \). Since \( x \) and \( y \) must be colored differently, suppose without loss of generality that \( x \) is colored in color 1 and \( y \) in color 2.  

Define \( U(x) \) and \( U(y) \) as neighborhoods of \( x \) and \( y \), respectively, such that no colors other than 1 and 2 appear in these neighborhoods. The intersection of \( U(y) \) with the circle of radius 1 centered at \( x \) must be entirely colored 2.  

Now, consider two points \( y_1 \) and \( y_2 \) on this intersection, positioned on opposite sides of \( y \). By Proposition \ref{l1.1}, the symmetric difference of the circles of radius 1 centered at \( y_1 \) and \( y_2 \) contains no points of color 2. Consequently, the intersection of this symmetric difference with \( U(x) \) must be entirely colored 1 (see Fig. \ref{fig:two_2c_pts}).  

However, this intersection contains interior points lying on a circle of radius 1 centered at \( y \), which contradicts Proposition \ref{l1.3}, since \(1 \in \overline{c}(y) \). 
\end{proof}

\begin{proposition} \label{l1.6}  
Let \( \omega = \omega(u)  \) be a unit circle centered at $u$, and let points \( x \) and \( y \) lie on \( \omega \) at a distance of 1 from each other. Suppose that at each of \( x \) and \( y \), one half-neighborhood in $\omega$ is colored 1 and the other half-neighborhood  is colored 2.   
Then, both inside and outside \( \omega \), the points \( x \) and \( y \) are adjacent to regions of colors other than 1 and 2.  
\end{proposition}  

\begin{proof}  
Assume that the half-neighborhood of \( x \) on the side of \( y \) along \( \omega \) is colored 1. Then:  
\begin{itemize}
\item The opposite half-neighborhood of \( x \) is colored 2.  
\item The half-neighborhood of \( y \) on the side of \( x \) is colored 1.  
\item The opposite half-neighborhood of \( y \) is colored 2.  
\end{itemize}

Now, consider the intersection of the sufficiently small neighborhood of \( x \) in \( \mathbb{R}^2 \) with the interior of the circle $\omega(y)$ of radius 1 centered at \( y \). Any point $x'$ in this intersection cannot be colored 2, since there is $y' \in \omega$, $c(y')=2$, $\|x'-y'\|=1$.  

Similarly, no point in the intersection of the neighborhood of \( x \) with the exterior of   $\omega(y)$ can be colored 1.  

Since the arc of  $\omega(y)$ is not a boundary, any neighborhood of \( x \)  must contain colors other than 1 and 2 in the both of the interior and exterior of $\omega$.  
\end{proof}

\section{Limitations on Chromaticity}  

In this section, we prove that in any valid 6-coloring of the plane, where region boundaries are Jordan curves that do not contain arcs of radius 1, no point can have chromaticity \( k \geq 4 \). Moreover, if the boundaries of the regions are polylines, it is always possible to recolor the neighborhoods of points lying on the boundary of more than three monochromatic regions in such a way that the coloring remains valid while ensuring that no point belongs to the closure of more than three monochromatic regions.  

\begin{proposition} \label{prop_gr_4}  
If Condition \ref{cond_curv1} holds for a 6-coloring, then no point can have chromaticity greater than 4.  
\end{proposition}  

\begin{proof}  
Suppose, for contradiction, that there exists a point \( x \) with multicolor \( \{1,2,3,4,5\} \). Consider the circle \( \omega(x) \) of radius 1 centered at \( x \). By Condition \ref{cond_curv1}, only a finite number of points on \( \omega(x) \) can be of a nontrivial multicolor, while all remaining points must be monochromatic and belong to the interior of their respective regions. Since the plane is colored with exactly six colors, these remaining points would all have to be assigned color 6. However, this contradicts the assumption that the coloring adheres to the given conditions.  

A similar contradiction arises when considering a point \( x \) of multicolor $\{1, \dots, 6\}$. By the same reasoning, the points of \( \omega(x) \) that do not lie on region boundaries cannot be assigned any valid color, leading to an impossibility. This establishes that no point can have chromaticity greater than 4.  
\end{proof}  

\begin{lemma} \label{t3}  
Under Condition \ref{cond_curv1}, no point in a valid 6-coloring can have chromaticity \( k \geq 4 \).  
\end{lemma}

\begin{proof}  
It remains to show that 4-chromatic points cannot exist. Suppose, for contradiction, that such a point \( u \) exists, and let \( \omega = \omega(u) \) be the circle of radius 1 centered at \( u \). Without loss of generality, assume \( u \) has multicolor \( \{3,4,5,6\} \). Then \( \omega \) must be divided into monochromatic arcs of colors \( 1 \) and \( 2 \).  

Let \( x_1 \) be a boundary point between monochromatic arcs of colors \( 1 \) and \( 2 \), and let \( x_2, x_3, x_4, x_5, x_6 \) be five additional points on \( \omega \) such that \( x_1, x_2, x_3, x_4, x_5, x_6 \) form a regular hexagon. By construction, each \( x_i \) is also a boundary point between the monochromatic arcs on \( \omega \).  

By Proposition \ref{l1.6}, both inside and outside \( \omega \), each \( x_i \) is adjacent to a region of a color different from \( 1 \) and \( 2 \). Furthermore, by Proposition \ref{l1.4}, the colors of these adjacent regions inside and outside \( \omega \) must be distinct.  

From Proposition \ref{prop_gr_4}, it follows that at each \( x_i \), exactly one region of a color from \( \{3,4,5,6\} \) is adjacent inside \( \omega \) and exactly one region of a different color from \( \{3,4,5,6\} \) is adjacent outside \( \omega \). Without loss of generality, assume that at \( x_1 \), the adjacent regions are of colors \( 3 \) and \( 4 \).  

By Proposition \ref{l1.3}, no point on the circle of radius 1 centered at \( x_2 \) can be a monochromatic point of color \( 1 \) or \( 2 \). Since arcs of circles cannot serve as boundaries of monochromatic regions, it follows that in the intersection of the circle centered at \( x_2 \) with radius 1 and the neighborhood of \( x_1 \), there exist monochromatic points of colors \( 3 \) and \( 4 \). Consequently, by Proposition \ref{l1.3}, colors \( 3 \) and \( 4 \) are not included in the multicolor of \( x_2 \), meaning that \( x_2 \) must be adjacent to regions of colors \( 5 \) and \( 6 \).  

Applying the same argument iteratively, we deduce that:
\begin{itemize}
\item \( x_3 \) and \( x_5 \) are adjacent to regions of colors \( 3 \) and \( 4 \),
\item \( x_4 \) and \( x_6 \) are adjacent to regions of colors \( 5 \) and \( 6 \).  
\end{itemize}

Thus, for each color in \( \{3,4,5,6\} \), there exist two points at an arc distance of \( \frac{2\pi}{3} \) along \( \omega \), where regions of that color are adjacent from the same side (either inside or outside \( \omega \)).  

Define \( \alpha_i \) as the half-plane containing \( x_i \) whose boundary is the line passing through \( u \) and perpendicular to \( ux_i \). If the region of a certain color is adjacent to \( x_i \) from the interior of \( \omega \), then the neighborhood of \( u \) colored in that color must lie in \( \alpha_i \); if it is adjacent from the exterior, then the corresponding neighborhood of \( u \) must lie in \( \alpha_{i+3} \).  

Since for each color there exist two points at an arc distance of \( \frac{2\pi}{3} \) where the corresponding regions are adjacent from the same side, it follows that the intersection of the neighborhood of \( u \) with each of the colors \( 3,4,5,6 \) is confined to a sector of size at most \( \frac{\pi}{3} \). This leads to a contradiction, as four such disjoint sectors cannot fit within a full neighborhood of \( u \).  

Thus, no point can have chromaticity 4, completing the proof.  
\end{proof}

\begin{figure}
    \centering
    \includegraphics[width=7cm]{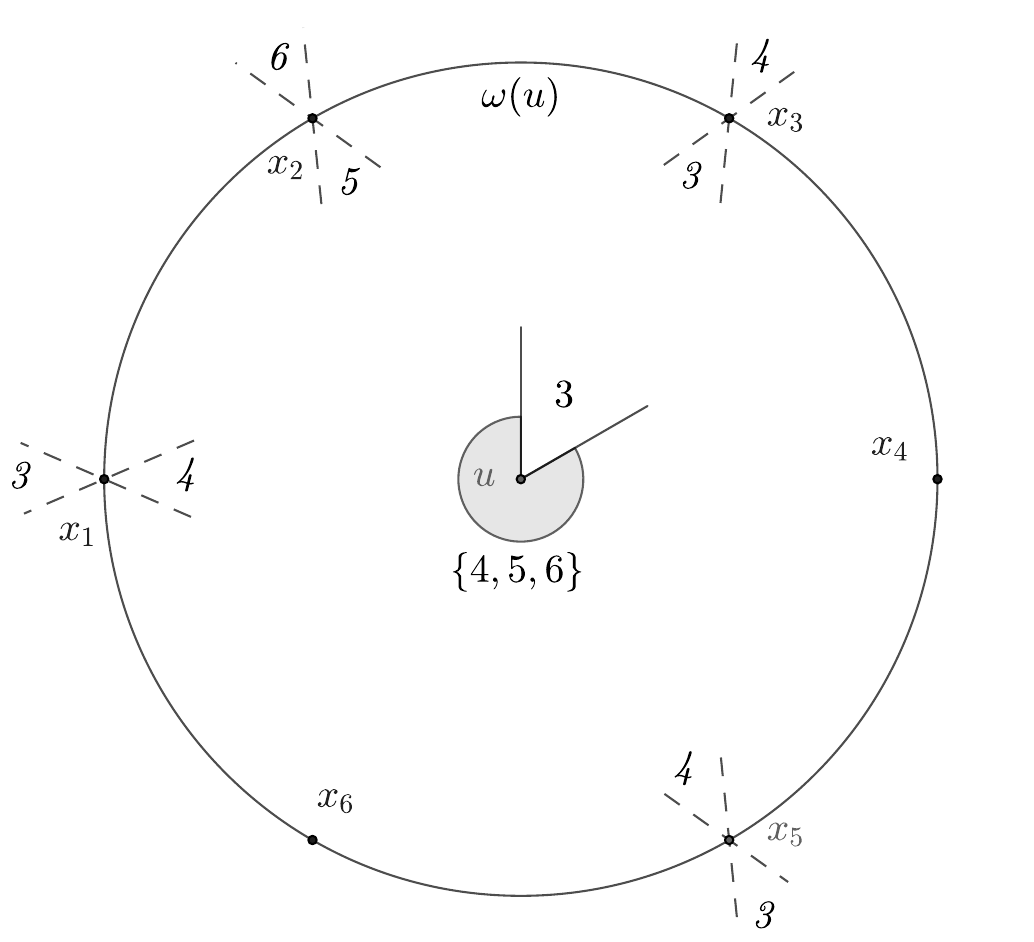} 
    \caption{Neighborhood of the unit circle centered at a 4-chromatic point}
    \label{fig_4col_point}
\end{figure}

\begin{proposition} \label{l6.1}
Let \( x, y, z \) be three non-collinear points in the plane, and suppose that all interior points of the segments \( zx \) and \( zy \) are colored \( 1 \). Then there exist points \( p \) and \( q \) on the segments \( zx \) and \( zy \), respectively, such that the coloring remains proper if all interior points of the triangle \( pqz \) are recolored with color \( 1 \).
\end{proposition}

\begin{proof} 
Consider a circle \( \omega(v) \) of unit radius passing through \( x \) and \( y \), with \( v \) and \( z \) lying on opposite sides of the line \( xy \). Draw a tangent \( l \) to \( \omega(v) \) that is parallel to \( xy \). Let \( p \) and \( q \) be the intersections of \( l \) with the segments \( xz \) and \( yz \), respectively.  

Now, consider any unit circle \( \omega(u) \) that intersects the interior of the triangle \( pqz \). Two cases arise:  

\begin{itemize}

        \item \( \omega(u) \) intersects at least one of the segments \( xz \) or \( yz \).
   In this case, by the properties of the given coloring, we must have \( c(u) \neq 1 \).  

        \item \( z \in \omega(u) \) and exactly one of the points \( x \) or \( y \) lies in the interior of \( \omega(u) \)
   By Proposition \ref{l1.4}, this implies that \( c(u) \neq 1 \).  
\end{itemize}

Since in both cases the color of \( u \) is not \( 1 \), we conclude that there is no obstruction to recoloring the interior of the triangle \( pqz \) with color \( 1 \). This completes the proof.  
\end{proof}

\begin{proposition} \label{t4}
If there exists a polygonal 6-coloring, then for any arbitrary \( r>0 \), there exists a polygonal coloring in 6 colors such that no point within a circle of radius \( r \) has more than three boundary segments meeting at it.
\end{proposition}

\begin{proof}
Fix an arbitrary circle of radius \( r \). If there exist trichromatic points within this circle where more than three boundary segments converge, we iteratively recolor the triangle whose existence is guaranteed by Proposition \ref{l6.1}. By applying this recoloring process, we ensure that no such points remain, while maintaining a valid 6-coloring.
\end{proof}

\section{Annulus Centered at the Trichromatic Point}

Let \( u \) be a trichromatic point, i.e., a common point of three monochromatic regions. Without loss of generality, assume that \( u \) has multicolor \( \{4, 5, 6\} \). We examine the coloring of the neighborhood of a circle of radius \( 1 \) centered at \( u \). We prove that if the boundaries between monochromatic regions do not contain an arc of a circle of radius \( 1 \), then in any neighborhood of this circle, there exists a closed curve colored in colors \( 1, 2, 3 \) enclosing \( u \). 

As before, denote by $\omega=\omega(u)$ the unit circle centered at $u$. Consider a circle $\omega_\delta(u)$ of sufficiently small radius \( \delta \) centered at \( u \). Define the points \( x_{45}, x_{46}, x_{56} \) as the intersections of $\omega_\delta(u)$ with the boundaries separating regions of colors  4 and 5, colors 4 and 6, colors 5 and 6, respectively. If a boundary intersects the circle multiple times, choose an arbitrary intersection point.

\begin{proposition} \label{l5}
Fix a ray \( l \) emanating from \( u \). Let \( \alpha \) be the second largest of the three angles between \( l \) and the vectors \( ux_{45}, ux_{46}, ux_{56} \). Define
\[
h_\delta = h_\delta(\alpha) = \delta \cos \alpha + \sqrt{1 - \delta^2 \sin^2\alpha}.
\]
Then the segment of this ray with endpoints at distances \( 1 \) and \( h_\delta \) from \( u \) is colored in colors \( 1, 2, 3 \).
\end{proposition}

\begin{proof}
A point located on the given ray at a distance \( h_\delta \) is at a distance of exactly \( 1 \) from the second farthest of the three points \( x_{45}, x_{46}, x_{56} \). Any point closer to \( u \) on this ray is at a distance less than \( 1 \) from at least two of these three points, while any point farther than \( 1 \) from \( u \) is at a distance greater than \( 1 \) from at least two of these points.

Thus, if \( h_\delta > 1 \), the points in the interval \( (1, h_\delta) \) are at a distance greater than \( 1 \) from \( u \) and less than \( 1 \) from two of \( x_{45}, x_{46}, x_{56} \). By Proposition \ref{l1.2}, these points cannot be assigned any of the colors \( 4,5,6 \), since \( u \) and at least  two of the three points lie on the boundary of these colors.

Similarly, if \( h_\delta < 1 \), the points in the interval \( (h_\delta, 1) \) are at a distance less than \( 1 \) from \( u \) and greater than \( 1 \) from two of the three points. By the same reasoning, they cannot be colored \( 4,5,6 \), ensuring that the interval is entirely colored in colors \( 1, 2, 3 \). 
\end{proof}

\begin{figure}
    \centering
    \includegraphics[width=8cm]{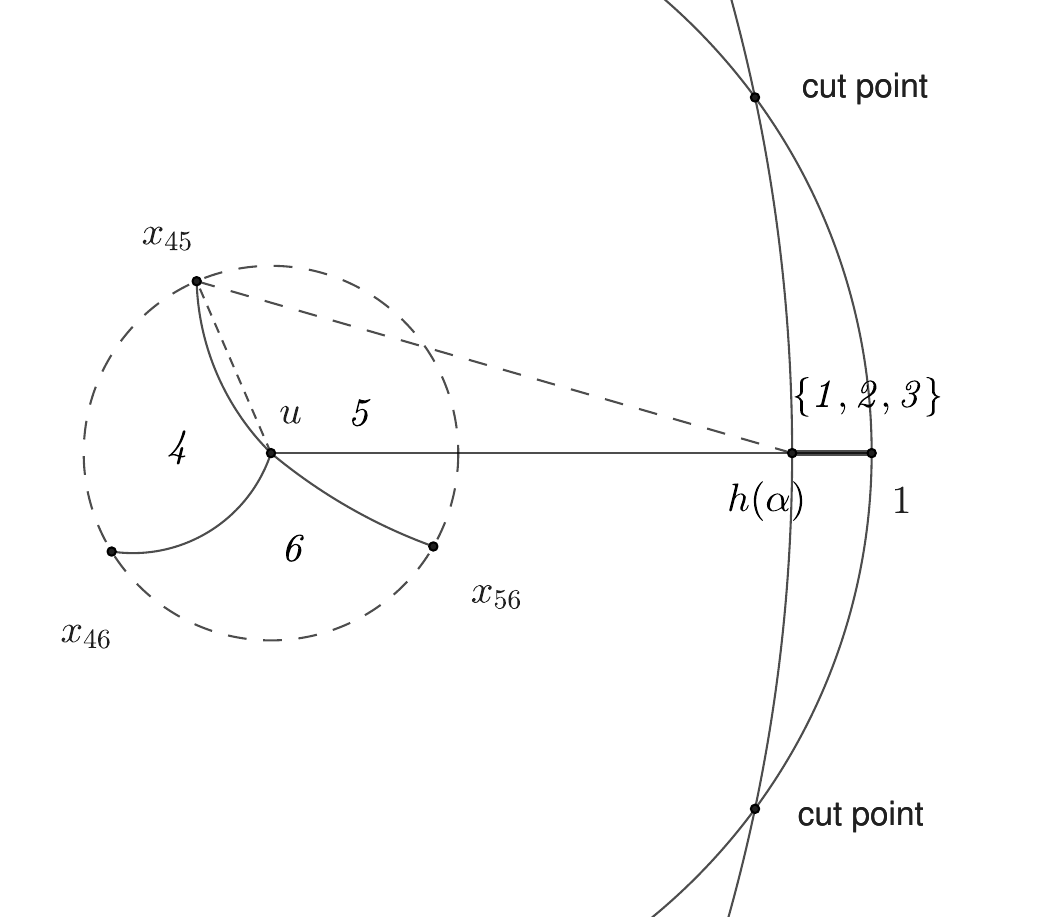}
    \label{fig_interval_3color}
    \caption{3-colored strip bordered with unit circle}
\end{figure}

Applying this lemma for all values of \( \alpha \), we conclude that in any neighborhood of a circle of radius \( 1 \), there exists an open strip colored in colors \( 1, 2, 3 \). This strip shrinks to a point (or possibly vanishes) in those directions where, in the notation of the lemma, \( h_\delta(v) = 1 \).

However, we can also adjust $\delta$. Define the sets
\[
A_{\delta} = \{ v \in \omega \mid h_\delta(v) < 1 \}, \quad
B_{\delta} = \{ v \in \omega \mid h_\delta(v) > 1 \}.
\]
Now, let
\[
A = \bigcup\limits_{0<\delta<\delta_1} A_\delta, \quad
B = \bigcup\limits_{0<\delta<\delta_1} B_\delta.
\]

\begin{proposition} \label{prop_interior}
Every point of \( \omega \) belongs to the interior of at least one of the sets \( A \) or \( B \), i.e., 
\[
\operatorname{Int} A \cup \operatorname{Int} B = \omega.
\]
\end{proposition}

\begin{proof}
Assume for contradiction that there exists a point \( v \) such that \( v \in \partial A \) and \( v \in \partial B \). This implies that for all \( 0<\delta<\delta_1 \), the unit circle centered at \( v \) passes through at least one of the points \( x_{45}(\delta) \), \( x_{46}(\delta) \), or \( x_{56}(\delta) \). 

Since Condition \ref{cond_loc_finite} (local finiteness) holds, it follows that at least one of the boundaries must contain an arc of a unit circle of positive length. However, this contradicts Condition \ref{cond_curv1}, which states that the boundaries do not contain such arcs. Hence, the assumption is false, proving that every point of \( \omega \) belongs to the interior of at least one of the sets \( A \) or \( B \). 
\end{proof}

In other words, for every trichromatic point \( u \), $\overline{c}(u)=\{a,b,c\}$, there exists a deformed annulus colored in the other three colors, whose closure contains the unit circle \( \omega = \omega(u) \). Next, we will choose one of many possible ways to show that the deviations of a piecewise smooth trichromatic curve from the circle can be made arbitrarily small.

\begin{defin}
Let \( t >0 \), and let \( \gamma \subset \omega \) be an arc of the circle. The annulus sector \( Q(u, \gamma, t) \) is a curved quadrilateral bounded by the arcs \( \gamma \), \( \gamma' = \{t x \mid x \in \gamma\} \), and the ray segments passing through the center \( u \) and the endpoints of \( \gamma \) and \( \gamma' \).
\end{defin}

From Proposition \ref{prop_interior}, we obtain the following corollary:

\begin{corollary} \label{prop_ann_sectors}
There exists \( \eta > 0 \) and finite families of arcs \( \Gamma_1=\{\gamma_{1,1}, \dots ,\gamma_{1,s}\} \) and \( \Gamma_2=\{\gamma_{2,1}, \dots ,\gamma_{2,s}\} \), where \( \gamma_{i,j} \subset \omega \), such that the following conditions hold:

\medskip
(a) The interiors of these arcs cover the entire circle:
\[
\bigcup\limits_{i=1,2,\; 1 \leq j \leq s} \operatorname{Int} \gamma_{i,j} = \omega.
\]

\medskip
(b) For each \( 1 \leq i \leq s \), the annulus sectors \( Q(u, \gamma_{1,i}, 1-\eta) \) and \( Q(u, \gamma_{2,i}, 1+\eta) \) are colored in colors \( 1,2,3 \).
\end{corollary}

\begin{figure}
    \centering
    \includegraphics[width=8cm]{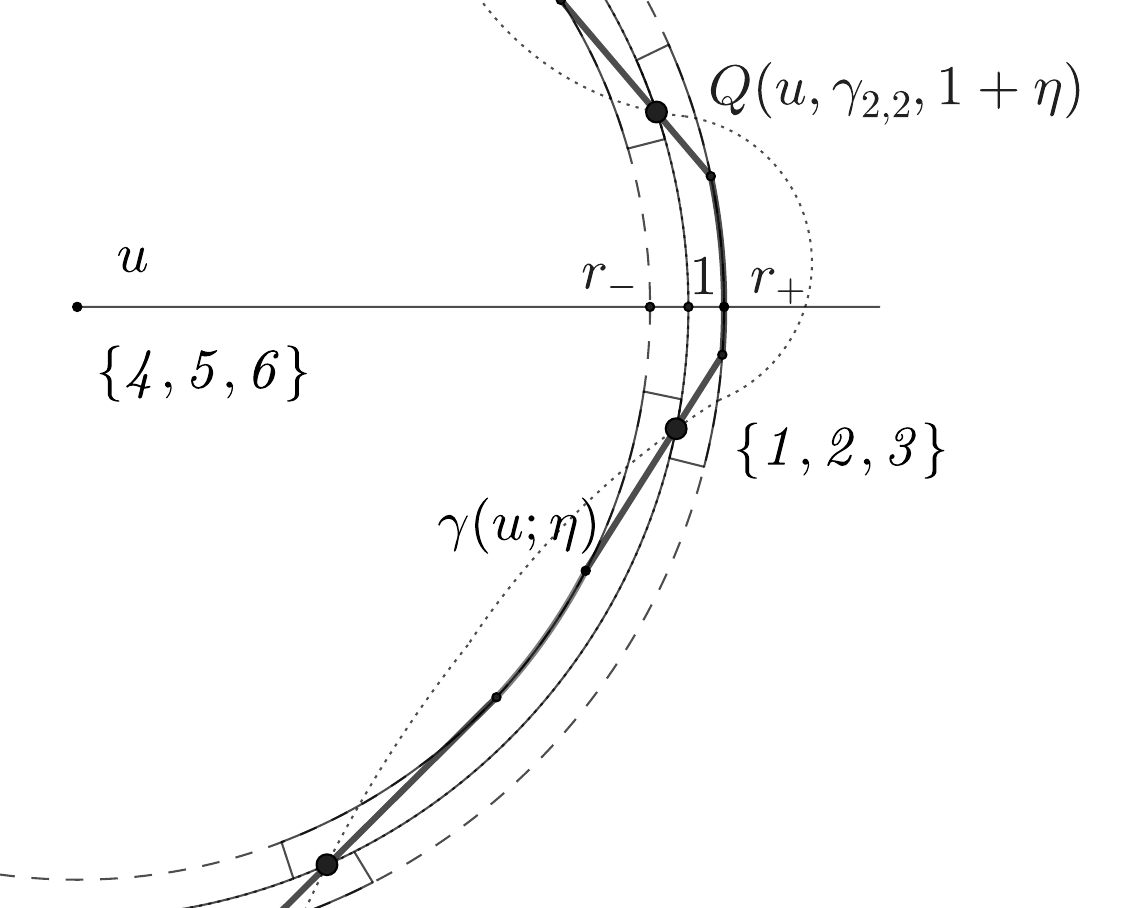}
    \caption{Construction of the piecewise-smooth curve in the interior of $C_1 \cup C_2 \cup C_3$}
    \label{fig:circle_strip}
\end{figure}

\begin{proposition} \label{cor_curve_exist}
By Corollary \ref{prop_ann_sectors}, for any \( \eta \), there exists a piecewise smooth closed curve \( \gamma(u; \eta) \) such that:
\begin{enumerate}
    \item The point \( u \) lies inside this curve.
    \item The curve is entirely contained within the union of regions colored \( 1, 2, 3 \).
    \item The curve remains within the \( \eta \)-neighborhood of \( \omega(u) \).
    \item At any point on the curve where the tangent is uniquely defined, the angle between the tangent and the direction to the point \( u \) lies within the interval \( (\pi/2 - \theta, \pi/2 + \theta) \), and it is possible to choose a curve such that \( \theta = o(\eta) \).
\end{enumerate}
\end{proposition}

\begin{proof}
Consider arcs lying outside the circle \( \omega \) as portions of a circle with radius \( r_+ \in (1, 1 + \eta) \), and arcs lying inside \( \omega \) as portions of a circle with radius \( r_- \in (1 - \eta, 1) \). Draw segments inclined at an angle \( \theta \), where \( 0 < \theta < \eta \), relative to the tangent at their points of intersection with these arcs (see Fig. \ref{fig:circle_strip}).
\end{proof}

\begin{remark}
We can assume that \( \theta = o(\eta^k) \) for any \( k \geq 1 \), though in future considerations, it will suffice to take \( k = 1 \).
\end{remark}

The following corollary is one of the key ingredients of the main result. Note that we cannot assert it unless the Condition \ref{cond_curv1} is satisfied and the breakpoints in rings may be unrecoverable.

\begin{corollary} \label{f3.2}
There do not exist two trichromatic points with disjoint multicolors such that the distance between them is less than \(2\) and the boundaries between each pair of colors in their respective multicolor sets extend from each point.
\end{corollary}

\begin{proof}
For sufficiently small \( \eta \), the curves constructed in Proposition \ref{cor_curve_exist} corresponding to these two trichromatic points must intersect. At the intersection points, none of the six colors can be assigned, leading to a contradiction.
\end{proof}

\subsection{Complementary Curves}

In this section, we present an alternative proof of the lemma on the coloring of curves marked by the endpoints of a unit segment and provide a reformulation of an auxiliary statement. In \cite{voronov2024chromatic}, this result was established for colorings with a forbidden distance interval.

\begin{defin}
Two curves \( \gamma_1, \gamma_2 \) are called \textit{complementary curves} with respect to a set of three colors \( M \) if all their points are interior points of the union of regions colored in \( M \), and there exist continuous parameterizations \( \gamma_1, \gamma_2: [0, 1] \to \mathbb{R}^2 \) such that for every \( t \in [0, 1] \), the condition
\[
\|\gamma_1(t) - \gamma_2(t)\| = 1
\]
is satisfied.
\end{defin}

\begin{defin}
Consider a line \( l_{1,2,3} = \{ x \mid x \in \mathbb{R} \} \), which does not belong to the given plane but is colored in three colors. The color of a point \( x \in l_{1,2,3} \) is determined by the formula
\[
c(x) = (\lfloor x + 1/2 \rfloor \operatorname{mod} 3) +1.
\]
Now, let \( f: \gamma([0,1]) \to l_{1,2,3} \) be a continuous mapping that preserves the coloring, i.e., \( c(f(x)) = c(x) \) for all points except possibly a finite number of interior points of the segment. Additionally, we require that \( f(\gamma(0)) \) and \( f(\gamma(1)) \) belong to \( \mathbb{Z} \). 

We define the \textit{index} of the curve \( \gamma \) as the integer
\[
\operatorname{Ind}(\gamma) = f(\gamma(1)) - f(\gamma(0)).
\]
\end{defin}

\begin{proposition}
    If  $\gamma_1 \sim \gamma_2$, then 
    \[
        |\operatorname{Ind} \gamma_1 - \operatorname{Ind} \gamma_2|\leq 1.
    \]
\end{proposition}

\begin{proof}
Let $\gamma(t_*)$ be a bichromatic point. We say that $t_*$ is an \emph{increasing color} point if, as $t$ increases through $t_*$, the color of the point $\gamma(t)$ changes cyclically as $1 \to 2$, $2 \to 3$, or $3 \to 1$.

We construct a piecewise linear function $f: \gamma([0,1]) \to l_{1,2,3}$ that preserves color. It suffices to define $f$ at bichromatic points and at the endpoints of the segment, i.e., on the set \( 0 = t_0, t_1, \dots, t_s = 1 \), as follows:

\[
f(\gamma(0)) = f(\gamma(t_0)) = c(\gamma(0)), \quad f(\gamma(1)) = f(\gamma(t_s)) = f(\gamma(t_{s-1}))
\]

\[
f(\gamma(t_{k+1})) = 
\begin{cases}
    f(\gamma(t_k)) + 1, & \text{if } t_{k+1} \text{ is an increasing color point,} \\
    f(\gamma(t_k)) - 1, & \text{otherwise,} 
\end{cases}
\]
for \( 1 \leq k \leq s-2 \).

Now, define the function \( F(t) = f_1(\gamma_1(t)) - f_2(\gamma_2(t)) \). Since both \( f_1 \) and \( f_2 \) are piecewise linear and preserve coloring, \( F \) is continuous.

We claim that \( F(t) \) cannot take integer values divisible by 3. Indeed, if \( F(t') \equiv 0 \pmod{3} \) at some \( t' \), then \( c(f_1(\gamma_1(t'))) = c(f_2(\gamma_2(t'))) \), which implies \( c(\gamma_1(t')) = c(\gamma_2(t')) \). However, this contradicts the assumption that \( \|\gamma_1(t') - \gamma_2(t')\| = 1 \), which ensures that \( \gamma_1(t') \) and \( \gamma_2(t') \) have different colors.

Since \( f_1(\gamma_1(0)) \), \( f_2(\gamma_2(0)) \), \( f_1(\gamma_1(1)) \), and \( f_2(\gamma_2(1)) \) are integers, we conclude that \( F(0) \) and \( F(1) \) are also integers. Given that \( F(t) \) cannot take values divisible by 3, the possible values at \( t=1 \) satisfy:

\begin{itemize}
    \item  If \( F(0) \equiv 1 \pmod{3} \), then \( F(1) \in \{F(0), F(0) + 1\} \).
\item If \( F(0) \equiv 2 \pmod{3} \), then \( F(1) \in \{F(0), F(0) - 1\} \).
\end{itemize}

In either case, we obtain:

\[
|\operatorname{Ind} \gamma_1 - \operatorname{Ind} \gamma_2| = 
\left| f_1(\gamma_1(1)) - f_1(\gamma_1(0)) - f_2(\gamma_2(1)) + f_2(\gamma_2(0)) \right| =
\]

\[
= |F(1) - F(0)| \leq 1.
\]
\end{proof}

\begin{corollary}
Let the curves $\gamma_1, \gamma_2, \dots, \gamma_s$ be colored in three colors and satisfy the following conditions:

\begin{itemize}
    \item[(a)] Their endpoints are not bichromatic.
    \item[(b)] The colors at the endpoints are different, i.e., $c(\gamma_i(0)) \neq c(\gamma_i(1))$ for all $i = 1,2, \dots, s$.
    \item[(c)] The curves are pairwise equivalent in sequence, i.e., $\gamma_i \sim \gamma_{i+1}$ for all $i = 1,2, \dots, s-1$.
\end{itemize}

Then the indices of the curves, $\operatorname{Ind} \gamma_1, \operatorname{Ind} \gamma_2, \dots, \operatorname{Ind} \gamma_s$, are either all positive or all negative. 
\label{ind_sequence}
\end{corollary}

\begin{proof}
By condition (b), each index $\operatorname{Ind} \gamma_i$ is nonzero. From the previous proposition, we know that for any two consecutive curves $\gamma_i$ and $\gamma_{i+1}$, their indices satisfy 
\[
|\operatorname{Ind} \gamma_i - \operatorname{Ind} \gamma_{i+1}| \leq 1.
\]
Since indices are nonzero, their signs must be the same for consecutive curves. By induction, this property extends to all $\gamma_1, \gamma_2, \dots, \gamma_s$.
\end{proof}

\subsection{Two trichromatic points with coinciding multicolors}

In this section, we introduce an auxiliary result concerning the intersection of a unit circle with a curve lying in the neighborhood of another unit circle. Speaking non-strictly, here we mean the following statement: if two smooth curves intersect at a positive angle, the number of intersection points will not change if the perturbation of one curve (or even both curves) is small enough.

Consider two points $u, v$ such that $1 < \|u - v\| < 2$. By Proposition \ref{cor_curve_exist},  construct a piecewise smooth closed curve $\gamma = \gamma(u, \eta)$, which lies within the $\eta$-neighborhood of $\omega(u)$. For a given point $q \in \omega(u)$, denote by $q'$ the intersection of the ray $uq$ with the curve $\gamma$. Define $\alpha(q)$ as the angle between the tangent to $\omega(u)$ at $q$ and the tangent to $\gamma$ at $q'$, provided that the tangent exists.

\begin{proposition}
Let $p$ be the intersection point of $\omega(u)$ and $\omega(v)$, and let $\alpha_*$ be the minimum angle at which the circles $\omega_r(u)$ of radius $r \in [1-\eta, 1+\eta]$ intersect the circle $\omega(v)$. Suppose that 
\[
\forall q \in \omega(u): \quad |\alpha(q)| < \alpha_*.
\]
Then in the $\eta$-neighborhood of $p$, there exists exactly one point where the curve $\gamma(u, \eta)$ intersects the circle $\omega(v)$.
\label{curve_intersect}
\end{proposition}

\begin{proof}
Consider a parameterization $\varphi: [0,1] \to \gamma$. Under the given conditions, the distance function $\rho(t) = \|\varphi(t) - v\|$ varies monotonically when $\varphi(t)$ remains within the $\eta$-neighborhood of $p$. This follows from the fact that the scalar product
\[
\left\langle \frac{d \varphi(t)}{dt}, \varphi(t) - v \right\rangle = \frac{1}{2} \frac{d}{dt} \|\varphi(t) - v\|^2
\]
preserves its sign throughout the  $\eta$-neighborhood of $p$. Consequently, the function $\rho(t)$ can cross the value $1$ at most once, implying that $\gamma(u, \eta)$ intersects $\omega(v)$ at exactly one point in this region.
\end{proof}

\begin{lemma} \label{t7}
Suppose that the Conditions \ref{cond_map}, \ref{cond_curv1} (``map-type'', ``forbidden arcs'') are satisfied for a proper 6-coloring. Then there do not exist two trichromatic points $u_1$, $u_2$ with same multicolor such that $1 <\|u_1-u_2\|<2$ and the boundaries between each pair of colors in their multicolor extend from each of them.
\end{lemma}

\begin{figure}
    \centering
    \includegraphics[width=12cm]{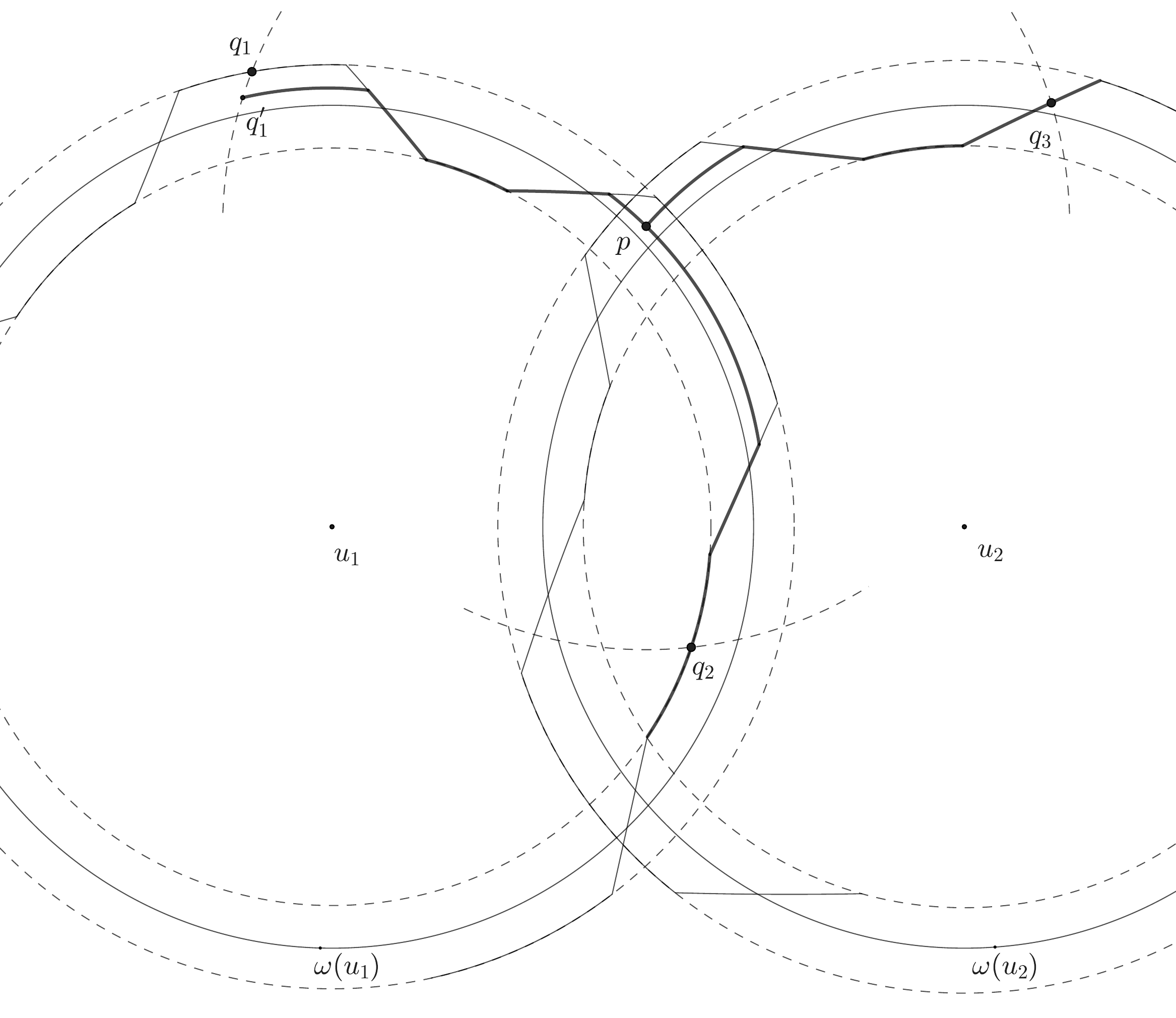}
    \caption{Construction of complementary piecewise-smooth trichromatic curves}
    \label{fig:tripod}
\end{figure}

\begin{proof}
Assume, for contradiction, that such trichromatic points $u_1, u_2$ exist. Without loss of generality, suppose both $u_1$ and $u_2$ have the multicolor $\{4,5,6\}$. 

In the neighborhoods of the circles $\omega(u_1)$ and $\omega(u_2)$, construct annuli $A(u_1)$ and $A(u_2)$, along with curves $\gamma(u_1, \eta_1) \subset A(u_1)$ and $\gamma(u_2, \eta_2) \subset A(u_2)$ that are colored with colors $1,2,3$. 

Next, choose a monochromatic point $p$ near the intersection point $p_0$. In the neighborhoods of the intersections of $\omega(p)$ with the curves $\gamma(u_1, \eta_1)$ and $\gamma(u_2, \eta_2)$, select monochromatic points $q'_1, q'_2, q'_3$ such that $q_1, q_3 \in A_1$ and $q_2 \in A_2$ (see Fig. \ref{fig:tripod}). This choice is possible because no arc of $\omega(p)$ of positive length can form a boundary. 

Now, consider the three curves $pq'_1$, $pq'_2$, and $pq'_3$ emanating from $p$ that satisfy the following condition:
\[
\forall x \in pq_i: \quad |\omega(x) \cap pq_j| = 1, \quad i,j \in \{1,2,3\}, \quad i < j.
\]

The uniqueness of the intersections is ensured by Proposition \ref{curve_intersect}. Under these conditions, it follows that the curves are complementary:
\[
q'_1 p \sim p q'_2, \quad q'_2 p \sim p q'_3, \quad q'_3 p \sim p q'_1.
\]

However, this leads to a contradiction in the coloring. Specifically, in a valid 3-coloring, we must have:
\[
\operatorname{Ind}(q'_1 p) = - \operatorname{Ind}(p q'_1).
\]
By Corollary \ref{ind_sequence}, the indices of these complementary curves must have the same sign. This contradiction implies that the original assumption was false.
\end{proof}

\section{Coloring a circle centered at a trichromatic point}.

In this section, we obtain the key ingredients to prove the following statement under the conditions of Theorems \ref{t2}, \ref{t1}. 

\begin{lemma} \label{t10}
Suppose that the Conditions \ref{cond_map}, \ref{cond_curv1}, \ref{cond_3col} (``map-type'', ``forbidden arcs'', ``no trichromatic vertices of degree 4+'') are satisfied for a proper 6-coloring. Then there are either two trichromatic points between $1$ and $2$ with the same multicolor, or two trichromatic points less than $2$ with non-intersecting multicolors.
\end{lemma}

Now, let us fix a coloring of a circle of radius $3$ that satisfies the conditions of Lemma \ref{t10}.

It is straightforward to verify that within a distance of at most $1$ from any point, there exists a trichromatic point. Fix such a trichromatic point $u$ located within a distance of at most $1$ from the center of the circle. Without loss of generality, assume that $u$ has the multicolor set $\{4, 5, 6\}$. 

Define $\omega$ as the circle of radius $1$ centered at $u$, and let $A(u)$ be the region described in Proposition \ref{prop_ann_sectors} that is colored with colors $\{1,2,3\}$ and whose closure contains $\omega$.

\subsection{Coloring the circle in three colors}

Next, we examine the coloring of the circle $\omega = \omega(u)$, which arises from the proper coloring of the plane under the conditions of Theorem \ref{t2} or Theorem \ref{t1}. Some properties of the coloring of a circle centered at a trichromatic point can be established independently from the plane’s coloring (Lemma \ref{l5.2}). In other cases, we will rely on previously established properties of the colorings of regions adjacent to $\omega(u)$.

By coloring the circle $\omega(u)$, we mean partitioning $\omega(u)$ into a finite number of arcs and assigning each arc one of three colors, ensuring that no two adjacent arcs share the same color. We define the interior points of these monochromatic arcs as \emph{monochromatic points}, which are assigned the color of the corresponding arc. The boundaries between adjacent arcs are called \emph{bichromatic points}, and they are considered to be colored in both colors of the adjoining arcs.

A coloring of $\omega$ is said to be \emph{proper} if:
\begin{itemize}
    \item No two monochromatic points at a distance of $1$ from each other share the same color.
    \item No two bichromatic points at a distance of $1$ from each other share the same pair of colors.
\end{itemize}

Let the center of $\omega(u)$ be a trichromatic point $u$ with multicolor $\{4,5,6\}$.

\begin{defin}
Let $c(x)$ denote the multicolor of a point $x \in \omega(u)$, as determined by the coloring of the entire plane. The \emph{pseudo-color} $c'(x)$ of a point $x$ is the set of colors given by $c'(x)=c(x) \cap \{1,2,3\}$.
\end{defin}

Let $A(u)$ be a region consisting of curved trapezoids colored in $\{1,2,3\}$ (see Corollary \ref{prop_ann_sectors}).

We first note a few obvious observations:

\begin{itemize}
    \item For all $x \in \omega(u)$, we have $c'(x) \neq \emptyset$, since $A(u)$ is colored with $\{1,2,3\}$ and $\omega(u)$ is contained in the closure of $A(u)$.
    
    \item A point has a two-element pseudo-color if and only if it lies on the boundary between colors passing through $A(u)$.
    
    \item No point can have the pseudo-color $\{1,2,3\}$; otherwise, it would be a trichromatic point with multicolor $\{1,2,3\}$ lying at distance $1$ from the trichromatic point $u$ with multicolor $\{4,5,6\}$, contradicting Corollary \ref{f3.2}.
    
    \item The circle $\omega$ intersects the color boundaries at a finite number of points, thus dividing $\omega$ into a finite number of monochromatic arcs. The endpoints of these arcs correspond to points with two-element pseudo-colors.
\end{itemize}

\begin{proposition}
No two points on $\omega$ at unit distance can have the same two-element pseudo-color.
\end{proposition}

\begin{proof}
Suppose, for contradiction, that two points $x, y \in \omega$ have the same pseudo-color $\{1,2\}$ and satisfy $\|x - y\| = 1$. By Proposition \ref{l1.6}, both inside and outside $\omega$, there exist regions adjacent to $y$ that are colored in colors other than $1$ and $2$. 

Since the intersection of the neighborhood of $y$ with either the interior or exterior of $\omega$ lies within a strip colored in $\{1,2,3\}$, at least one of the regions adjacent to $y$ must be colored $3$. This contradicts our assumption that no point on $\omega$ has all three colors $\{1,2,3\}$ in its neighborhood. Hence, no two points on $\omega$ at unit distance can have the same two-element pseudo-color.
\end{proof}

Thus, we arrive at the following conclusion:

\begin{proposition}\label{l5.1}
The coloring of the points of $\omega$ according to their pseudo-colors constitutes a proper coloring of $\omega$.
\end{proposition} 

Next, we analyze the possible colorings of $\omega$ independently of the plane.

\begin{defin}
A proper coloring of $\omega$ using three colors is called \emph{cyclic} if no monochromatic arc is bordered on both sides by arcs of the same color. That is, as we traverse the circle, the colors of the arcs follow a repeating sequence of either $1,2,3,1,2,3,\dots$ or $1,3,2,1,3,2,\dots$.
\end{defin}

\begin{lemma} \label{l5.2}
If $\omega$ is properly colored with three colors, then it is possible to recolor some of the monochromatic arcs of $\omega$ so that the coloring remains proper and becomes cyclic. This recoloring does not create new bichromatic points, and the previous bichromatic points either become interior points of monochromatic arcs or remain unchanged.
\end{lemma}

\begin{proof}

Consider the shortest monochromatic arc that is bordered on both sides by arcs of the same color. Such an arc must exist since the coloring is not cyclic. If multiple such arcs exist, choose any one of them. Without loss of generality, assume this arc is colored $1$, and its neighboring arcs are colored $2$. Let $x_1$ and $y_1$ be the endpoints of this arc in the clockwise direction.

If no point of the arc $x_1y_1$ is at distance $1$ from a point of color $2$, then we can simply recolor the entire arc $x_1y_1$ with color $2$. Otherwise, let $z_1$ be a point on $x_1y_1$ that is at distance $1$ from a point $w_1$ of color $2$. Without loss of generality, assume that $w_1$ is located clockwise from $z_1$. 

Now, define points $x_2, \dots, x_6$ on $\omega$ such that $x_1x_2x_3x_4x_5x_6$ forms a regular hexagon whose vertices are listed in clockwise order. Similarly, define $y_1, \dots, y_6$. Since $w_1$ lies on arc $x_2y_2$ and is colored $2$, the arc of color $2$ containing $w_1$ must be nested within the arc $x_2y_2$. This follows from the fact that the arc $x_1y_1$ is bordered on both sides by arcs of color $2$ and also on both sides by arcs of color $3$. Otherwise, one of its endpoints would have multicolor $\{1,2\}$ and be at distance $1$ either from an interior point of $x_1y_1$  or from the other endpoint of $x_1y_1$, which also has multicolor $\{1,2\}$ (contradicting the definition of a proper coloring).

Thus, the arc of color $2$ containing $w_1$ must be at least as short as $x_1y_1$ and therefore coincides with the arc $x_2y_2$. If arc $x_2y_2$ cannot be recolored with color $3$ while preserving the proper coloring, then there must exist a point $z_2$ on it that is at distance $1$ from a point $w_2$ of color $3$. Since arc $x_1y_1$ does not contain any points of color $3$, $w_2$ must lie on arc $x_3y_3$. 

Repeating this reasoning, if no arc can be recolored in the color of its neighboring arcs, then we obtain the following cyclic structure:
\begin{itemize}
\item Arc $x_3y_3$ is colored $3$ and is bordered by arcs of color $1$.
\item Arc $x_4y_4$ is colored $1$ and is bordered by arcs of color $2$.
\item Arc $x_5y_5$ is colored $2$ and is bordered by arcs of color $3$.
\item Arc $x_6y_6$ is colored $3$ and is bordered by arcs of color $1$.
\end{itemize}

In this case, we perform a cyclic recoloring as follows:
\begin{itemize}
\item Recolor arcs $x_1y_1$ and $x_4y_4$ with color $2$.
\item Recolor arcs $x_2y_2$ and $x_5y_5$ with color $3$.
\item Recolor arcs $x_3y_3$ and $x_6y_6$ with color $1$.
\end{itemize}

After this recoloring, the points $x_i$ and $y_i$ become interior points of the monochromatic arcs and are no longer bichromatic. 

Since the number of monochromatic arcs in $\omega$ is finite, by repeating these steps we will get a cyclic coloring in a finite number of steps.
\end{proof}

Next, we must show separately that the case with six bichromatic points on the circle $\omega(u)$ centered at the trichromatic point $u$ is impossible.

As before, let $A(u)$ be a region consisting of curved trapezoids colored with colors $1,2,3$ (see Corollary \ref{prop_ann_sectors}).

\begin{proposition} \label{l5.4}
There cannot exist three points $x, y, z$ on $\omega$ such that:
\begin{itemize}
    \item The point $y$ is at distance $1$ from both $x$ and $z$.
    \item The pseudo-color of $y$ is $\{1,2\}$
    \item the half-neighborhood of $x$ on the $y$-side colored in color $1$ and the half-neighborhood of $z$ on the $y$-side colored in color $2$.
\end{itemize}
\end{proposition}

\begin{proof} 
Consider the intersection of the neighborhood of $y$ with the exterior of the circle centered at $x$ with radius $1$. This region contains no points of color $1$. Similarly, the intersection of the neighborhood of $y$ with the exterior of the circle centered at $z$ with radius $1$ contains no points of color $2$. 

By assumption, the intersection of any neighborhood of $y$ with the exterior of $\omega$ lies in $A(u)$, which is colored with colors $1,2,3$. Since the neighborhood of $y$ outside $\omega$ cannot contain only colors $1$ and $2$, it must contain a point of color $3$. However, this contradicts the assumption that no point on $\omega$ has a multicolor of $\{1,2,3\}$.

Thus, such a configuration of points $x, y, z$ on $\omega$ is impossible.
\end{proof}

\begin{proposition} \label{l5.5}
On $\omega$, there do not exist six points $x_1, \ldots, x_6$ with two-element pseudo-colors such that these points are the vertices of a regular hexagon, where:
\begin{itemize}
    \item At $x_1$ and $x_4$, the semicircles on $\omega$ are colored $1$ and $2$.
    \item At $x_2$ and $x_5$, the semicircles on $\omega$ are colored $2$ and $3$.
    \item At $x_3$ and $x_6$, the semicircles on $\omega$ are colored $3$ and $1$.
    \item At each $x_i$, the semicircle on the side of $x_{i-1}$ is colored in the first color of the pair, and the semicircle on the side of $x_{i+1}$ is colored in the second.
\end{itemize}
\end{proposition}

\begin{proof}
Suppose that such a configuration exists. By the construction of $A(u)$, the neighborhood of $u$ is divided into three sectors by the segments $u x_{45}$, $u x_{46}$, and $u x_{56}$. Each of these angles $\alpha_{45}$, $\alpha_{46}$, $\alpha_{56}$ has an associated arc centered at the intersection of its bisector with $\omega$. The neighborhood of such an arc's exterior intersection with $\omega$ lies entirely in $A(u)$. The length of this arc is given by $\pi - \alpha_{ij} - \varepsilon$, where  $\varepsilon$ depends on $\delta$ (the radius of the circumcircle of the triangle $x_{45}x_{46}x_{56}$) and tends to zero as $\delta \to 0$. 

If the three angles are not equal, then the smallest of them corresponds to an arc of length greater than $\frac{\pi}{3}$. This means that one of the vertices $x_1, \ldots, x_6$ must be contained within such an arc, contradicting Proposition~\ref{l5.4}. 

If all three angles are equal, then the points $x_1, \ldots, x_6$ must be the endpoints of arcs of length $\frac{\pi}{3}$ centered at the extensions of the region boundaries emanating from $u$ (otherwise, we would again reach a contradiction with Proposition~\ref{l5.4}).

Now, consider the vertices lying on the arc centered at the intersection of $\omega$ with the extension of $u x_{56}$. Without loss of generality, suppose these are the vertices $x_1$ and $x_2$. Define $U(x_1)$ and $U(x_2)$ as the intersections of the neighborhoods of $x_1$ and $x_2$, respectively, with the exterior of two circles of radius $1$ centered at adjacent vertices of the hexagon. 

Using the same reasoning as in Proposition~\ref{l5.4}, we deduce that $U(x_1)$ and $U(x_2)$ must be colored with colors $4,5,6$. However, the circles centered at points in $U(x_1)$ and $U(x_2)$ intersect the regions colored $5$ and $6$ in the neighborhood of $u$, implying that they must be colored $4$. 

Finally, similar to the proof of Proposition~\ref{l5.1}, we conclude that the boundary between colors $4$ and $2$ in the neighborhood of $x_1$ must form an arc of a circle of radius $1$. A contradiction with Condition~\ref{cond_curv1}. 
\end{proof}

The following two propositions and Lemma \ref{l5.6} in this setting are proved in the same way as the analogous statements 3.2 and 3.3 in \cite{voronov2024chromatic}. For completeness, the proofs are given in the Appendix.

\begin{proposition}[3.2, \cite{voronov2024chromatic}] \label{circle_triple_greater_than_1}
If there are at least 9 bichromatic points in a cyclic regular coloring of a circle, then 
for any pair of colors $\{1,2\}$, $\{2,3\}$, $\{1,3\}$, $\{1,3\}$ there are three bichromatic points with such multicolor whose pairwise distances between them are greater than $1$.
\end{proposition}

\begin{proposition}[3.3, \cite{voronov2024chromatic}]
If a closed circle of unit radius centered at point $u$ is 6-colored under the conditions of Theorem \ref{t2}, and the circle $\omega(u)$ contains at least 9 bichromatic points, then for each pair of colors $\{1,2\}$, $\{2,3\}$, $\{1,3\}$, $\{1,3\}$ there are three bichromatic points with such a multicolor for which

(a) pairwise distances are greater than $1$, 

(b) the bichromatic boundaries extending from their point on the circle end at a tricolor point inside the circle.
\label{circle_triple_borders}
\end{proposition}

\begin{lemma} \label{l5.6}
For any pair of colors from $\{1,2\}$, $\{2,3\}$, and $\{3,1\}$ on $\omega$, there exist three points with that pseudo-color that are more than $1$ apart from each other.
\end{lemma}

\begin{proof}
Consider the coloring of $\omega$ defined by pseudo-colors. By applying the recoloring process from Lemma \ref{l5.2}, we obtain a cyclic coloring. By Proposition \ref{l5.5}, this coloring cannot consist of exactly six monochromatic arcs of length $\frac{\pi}{3}$. Consequently, by Proposition \ref{circle_triple_borders}, there exist three bichromatic points with the required multicolor at pairwise distances greater than $1$. 

Since the recoloring in Lemma \ref{l5.2} does not introduce new bichromatic points, these points must have been present in the original coloring.
\end{proof}

\section{Proof of the Main Result}

We now combine the previously established results to complete the proofs of Theorems \ref{t2} and \ref{t1}.

\textbf{Proof of Theorem \ref{t2}.} 
Assume that there exists a Jordan-measurable coloring of the plane with six colors such that any edge intersects any circle of radius $1$ at only a finite number of points, and no trichromatic vertex has degree greater than three. 

By Lemma \ref{t3}, there are no four-chromatic points, allowing us to apply Lemma \ref{t10} to any circle of radius $3$ inside this region. As a result, we obtain either:
\begin{itemize}
    \item Two trichromatic points corresponding to colors $1$ and $2$ with the same multicolor assignment, contradicting Lemma~\ref{t7}, or
    \item Two trichromatic points at a distance of less than $2$ with disjoint multicolor sets, contradicting Corollary~\ref{f3.2}.
\end{itemize}
Thus, such a coloring cannot exist. \qed

\textbf{Proof of Theorem \ref{t1}.} 
Assume that there exists a polygonal coloring of the plane with six colors. By Lemma \ref{t4}, for any sufficiently large radius $r$, there exists a six-color polygonal coloring in which no more than three boundary segments meet at any point within a circle of radius $r$. 

By Lemma \ref{t3}, this coloring contains no four-chromatic points. Applying Lemma \ref{t10} to any circle of radius $3$ inside this region, we obtain either:
\begin{itemize}
    \item Two trichromatic points corresponding to colors $1$ and $2$ with the same multicolor assignment, contradicting Lemma~\ref{t7}, or
    \item Two trichromatic points at a distance of less than $2$ with disjoint multicolor sets, contradicting Corollary~\ref{f3.2}.
\end{itemize}
Hence, such a coloring is not possible. \qed

 \section{Conclusion}

Note that all reasoning remains valid if we consider the interior of a circle of radius~3. 

 \begin{question}
What is the minimum area of a convex region for which we can assert, under the conditions of Theorem \ref{t2}, that a 6-coloring does not exist?
 \end{question}

 Certainly, the conditions of Theorem \ref{t2} leave some frustration.

 \begin{question}
 Is it true that 7 colors are needed to color a locally finite map in the plane if the boundaries can be arcs of unit circles and if trichromatic vertices of degree greater than 3 are allowed?
 \end{question}

The conditions of Theorems 1, 2 allow us to handle some exceptional cases and to carry out the rest of the reasoning in the same way as it was done for the forbidden distance interval. In the above formulations, all ``features'' form a locally finite set of points in the plane. If we exclude the neighborhoods of vertices and consider a bounded region, we can find a forbidden distance interval. It may seem strange that the existence of arcs of unit curvature (or of a countable number of intersections of a boundary with such an arc) is a problem of principle. But in this case, ``pecularity'' has Hausdorff dimension 1, not 0. In other words, in the absence of such arcs, the area of regions $S = S(\varepsilon)$, after excluding which we obtain a proper coloring with forbidden distance interval $[1-\varepsilon, 1+\varepsilon]$, at $\varepsilon\to 0$ has asymptotics
$$S(\varepsilon)=\Omega(\varepsilon^2),$$
and if arcs of unit curvature are present, then
$$S(\varepsilon)=\Omega(\varepsilon^1).$$

Note that if any upper bound $\mu(C_1 \cup \dots \cup C_6)\leq \mu^*_6 <1$ on the fraction of the plane that can be correctly colored in 6 colors with forbidden interval $[1 - \varepsilon, 1+ \varepsilon]$, $\varepsilon>0$ were established, it would follow definitively that $\chi_{map}(\mathbb{R}^2)=7.$ Unfortunately, we only know that there exists a 6-coloring in which the $1/6997$ part of the  plane remains uncolored \cite{parts2020what}.

\newpage

\section{Appendix. Proof of Article{\cite{voronov2024chromatic}}}.

\textbf{Proof of the Proposition\ref{circle_triple_greater_than_1}}

Here we consider the coloring of the circle separately from the coloring of the rest of the plane.

Assume the opposite. Then after switching to cyclic coloring there are two points $p_1$, $p_2$ of multicolor $\{1,2\}$ on the circle, the distance between which is less than 1. Let $p_2$ be located clockwise from $p_1$.

We prove that then there is another point of multicolor $\{1,2\}$ at a distance less than 1 from $p_2$, clockwise from $p_2$. Indeed, the arc $p_1 p_2$ contains at least one pair of adjacent bichromatic points $q_1$, $q_2$ of multicolors $\{2,3\}$, $\{3,1\}$, respectively. Then the points $q'_1$, $q'_2$ obtained by a clockwise rotation by $\pi/3$ have colors 1, 2, and between them we find the point $p_3$ of multicolor $\{1,2\}$.

This means that on each arc of length $\pi/3$ there is a multicolor point $\{1,2\}$. Let's divide the circle into 6 such arcs and consider three arcs with pairwise distances between them equal to 1. On each of them we will find a multicolor point $\{1,2\}$, from which we get the required point.

\textbf{Proof of Proposition{\ref{circle_triple_borders}}}.

Here we consider the coloring of a circle and its insides separately from the rest of the plane.

Consider all bichromatic multicolor points $\{1,2\}$ on a circle. The bichromatic boundaries containing them can either end at a trichromatic point inside the circle, or connect two bichromatic points. For any boundary $p_1 p_2$ connecting two bichromatic points $p_1$, $p_2$ on the circle, it is possible to completely recolor the region of diameter less than 1 that $p_1 p_2$ separates from the circle into one of the colors, after which the points $p_1, p_2$ become monochromatic. The coloring of the circle and the interior will remain correct. 

Let's perform all such operations. After that, let us consider the coloring of the circle separately from the interior, pass to cyclic coloring and use the Proposition \ref{circle_triple_greater_than_1}. Then for the found triple points, the boundaries end inside the circle and the pairwise distances between the points are greater than 1. For multicolors $\{2,3\}$, $\{1,3\}$ the reasoning is similar. \qed

\begin{figure}
    \centering
    \includegraphics[width=6cm]{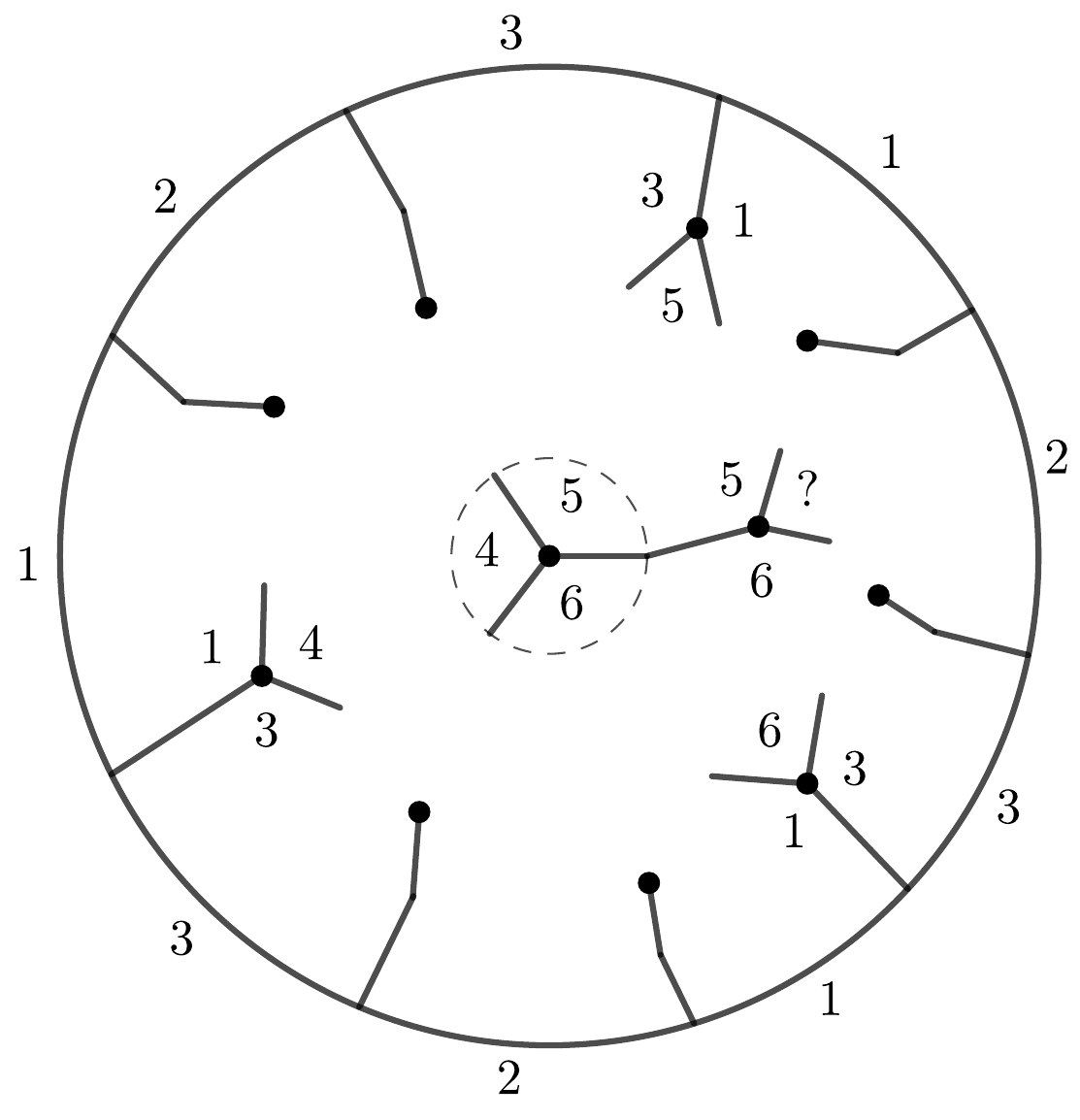}
    \caption{10+1 trichromatic points inside the unit circle}
    \label{fig:pts10}
\end{figure}

\textbf{Proof of Lemma \ref{t10}}.

Let the center of the circle have multicolor $\{4,5,6\}$. Let $B$ denote the set of 9 bichromatic points on the circle from the formulation of the theorem. Then $B$ contains 3 points of each of the multicolors $\{1,2\}$, $\{2,3\}$, $\{3,1\}$. Under the above conditions, any boundary that intersects the circle at a point from set $B$ must end at a tricolor point inside the circle. Otherwise, the boundary connects two bichromatic points with the same multicolor on the circle, but this is impossible because the distance between such points is greater than 1. Thus, in addition to $u$, we have 9 trichromatic points $p_1$, ... , $p_9$ of multicolors $\{1,2,i\}$, $\{1,3,i\}$, $\{2,3,i\}$, $i=4,5,6$. If the multicolors of two points $p_j$ coincide, we immediately get the desired one. Otherwise, inside the circle there are $\frac{1}{2}\binom{6}{3}=10$ trichromatic points whose multicolors are pairwise different and do not complement each other up to $\{1,2 \dots , 6\}$. It remains to notice that on the boundary of the union of three regions of colors $1,2,3$ whose vertex is $u$, there is at least one more trichromatic point (Fig. \ref{fig:pts10}). \qed

\bibliographystyle{abbrv}
\bibliography{main.bib}

\begin{thebibliography}{10}

\bibitem{warsaw}
J.~Chybowska-Sok{\'o}{\l}, K.~Junosza-Szaniawski, and K.~W{\k{e}}sek.
\newblock Coloring distance graphs on the plane.
\newblock {\em Discrete Mathematics}, 346(7):113441, 2023.

\bibitem{coulson2004chromatic}
D.~Coulson.
\newblock On the chromatic number of plane tilings.
\newblock {\em Journal of the Australian Mathematical Society}, 77(2):191--196, 2004.

\bibitem{deGrey}
A.~D. N.~J. de~Grey.
\newblock The chromatic number of the plane is at least 5.
\newblock {\em Geombinatorics}, 28(1):18--31, 2018.

\bibitem{exoo1}
G.~Exoo.
\newblock $\varepsilon$-unit distance graphs.
\newblock {\em Discrete and Computational Geometry}, 33(1):117--123, June 2004.

\bibitem{exoo2020chromatic}
G.~Exoo and D.~Ismailescu.
\newblock The chromatic number of the plane is at least 5: a new proof.
\newblock {\em Discrete \& Computational Geometry}, 64(1):216--226, 2020.

\bibitem{falconer1981realization}
K.~J. Falconer.
\newblock The realization of distances in measurable subsets covering $\mathbb{R}^n$.
\newblock {\em Journal of Combinatorial Theory, Series A}, 31(2):184--189, 1981.

\bibitem{manta2021triangle}
M.~N. Manta.
\newblock Triangle colorings require at least seven colors.
\newblock {\em Discrete Mathematics}, 344(7):112411, 2021.

\bibitem{parts2020what}
J.~Parts.
\newblock What percent of the plane can be properly 5- and 6-colored?
\newblock {\em Geombinatorics}, 30(1):25--39, 2020.

\bibitem{Soifer}
A.~Soifer.
\newblock {\em The mathematical coloring book: Mathematics of coloring and the colorful life of its creators}.
\newblock Springer Science \& Business Media, 2008.

\bibitem{townsend2005colouring}
S.~P. Townsend.
\newblock Colouring the plane with no monochrome unit.
\newblock {\em Geombinatorics}, 14(4):184--193, 2005.

\bibitem{voronov2024chromatic}
V.~Voronov.
\newblock The chromatic number of the plane with an interval of forbidden distances is at least 7, 2024.

\bibitem{woodall1973distances}
D.~R. Woodall.
\newblock Distances realized by sets covering the plane.
\newblock {\em Journal of Combinatorial Theory, Series A}, 14(2):187--200, 1973.

\end{thebibliography}

\end{document}